\documentclass[a4paper,11pt]{amsart}
\usepackage[T1]{fontenc}
\usepackage{url}
\usepackage{xcolor}
\usepackage{amsmath,amssymb}
\usepackage{stmaryrd}
\usepackage{nicefrac}
\usepackage[curve,matrix,arrow,pdf]{xy}
\usepackage[abbrev,non-sorted-cites]{amsrefs}
\BibSpec{article}{%
+{}{\PrintAuthors} {author}
+{,}{ \textit} {title}
+{,}{ } {journal}
+{}{ \textbf} {volume}
+{}{ \parenthesize} {date}
+{,}{ } {pages}
+{,}{ } {note}
+{.}{} {transition}
}
\BibSpec{collection.article}{%
    +{}  {\PrintAuthors}                {author}
    +{,} { \textit}                     {title}
    +{.} { }                            {part}
    +{:} { \textit}                     {subtitle}
    +{,} { \PrintContributions}         {contribution}
    +{,} { \PrintConference}            {conference}
    +{}  {\PrintBook}                   {book}
    +{, in} { \textit}                  {booktitle}
    +{,} { \PrintEditorsB}              {editor}
    +{, Vol.} { }                       {volume}
    +{,} { \PrintDateB}                 {date}
    +{,} { pp.~}                        {pages}
    +{,} { }                            {status}
    +{,} { available at \eprint}        {eprint}
    +{}  { \PrintTranslation}           {translation}
    +{;} { \PrintReprint}               {reprint}
    +{.} { }                            {note}
    +{,} { }                            {publisher}
    +{.} {}                             {transition}
}
\BibSpec{book}{%
    +{}  {\PrintPrimary}                {transition}
    +{.} { \PrintDate}                  {date}
    +{.} { \textit}                     {title}
    +{.} { }                            {part}
    +{:} { \textit}                     {subtitle}
    +{,} { \PrintEdition}               {edition}
    +{}  { \PrintEditorsB}              {editor}
    +{,} { \PrintTranslatorsC}          {translator}
    +{,} { \PrintContributions}         {contribution}
    +{,} { }                            {series}
    +{,} { \voltext}                    {volume}
    +{,} { }                            {publisher}
    +{,} { }                            {organization}
    +{,} { }                            {address}
    +{,} { }                            {status}
    +{}  { \parenthesize}               {language}
    +{}  { \PrintTranslation}           {translation}
    +{;} { \PrintReprint}               {reprint}
    +{.} { }                            {note}
    +{.} {}                             {transition}
    +{.} { }                            {url}
    +{}  {\SentenceSpace \PrintReviews} {review}
}
\BibSpec{thesis}{%
    +{}  {\PrintAuthors}                {author}
    +{.} { \PrintDate}                  {date}
    +{.} { \PrintDate}                  {year}
    +{.} { \textit}                     {title}
    +{:} { \textit}                     {subtitle}
    +{,} { \PrintThesisType}            {type}
    +{,} { }                            {school}
    +{,} { }                            {address}
    +{,} { \eprint}                     {eprint}
    +{,} { }                            {status}
    +{}  { \parenthesize}               {language}
    +{}  { \PrintTranslation}           {translation}
    +{;} { \PrintReprint}               {reprint}
    +{.} { }                            {note}
    +{.} {}                             {transition}
}
\usepackage[verbose]{hyperref}

\hypersetup{colorlinks=false,allbordercolors=blue,pdfborderstyle={/S/U/W 1}}

\renewcommand{\ge}{\geqslant}
\renewcommand{\le}{\leqslant}
\def\Cl#1{\ensuremath{\mathcal{#1}}}

\def\pv#1{\ensuremath{\mathsf{#1}}}

\def\Om#1#2{\ensuremath{\overline{\Omega}_{#1}\pv{#2}}}

\def\oms#1#2{\ensuremath{\Omega_{#1}^\sigma\pv{#2}}}

\def\omc#1#2{\ensuremath{\Omega_{#1}^\kappa\pv{#2}}}
\DeclareMathOperator{\Pol}{Pol}
\DeclareMathOperator{\BPol}{BPol}
\DeclareMathOperator{\UPol}{UPol}
\DeclareMathOperator{\BUPol}{BUPol}

\let\op=\llbracket
\let\cl=\rrbracket

\def\embedsin{\hookrightarrow}
\newcommand\malcev{%
\mathbin{\raise.5pt\hbox{$\bigcirc$\kern-8.5pt\raise.5pt\hbox{\scriptsize$m$}\,}}}

\begin{document}

\setcounter{page}{1}

\markboth{Almeida}{Pseudovarieties of semigroups}

\title{Pseudovarieties of semigroups}

\author{Jorge Almeida}

\address{CMUP, Faculdade de Ci\^encias, Universidade do Porto, Rua do
  Campo Alegre s/n, 4169-007
  Porto, Portugal}
\email{jalmeida@fc.up.pt}

\keywords{finite semigroup, pseudovariety, regular language,
  profinite semigroup}

\subjclass[2020]{20M07, 20M35}
\begin{abstract}
  The most developed aspect of the theory of finite semigroups is
  their classification in pseudovarieties. The main motivation for
  investigating such entities comes from their connection with the
  classification of regular languages via Eilenberg's correspondence.
  This connection prompted the study of various natural operators on
  pseudovarieties and led to several important questions, both
  algebraic and algorithmic. The most important of these questions is
  decidability: given a finite semigroup is there an algorithm that
  tests whether it belongs to the pseudovariety? Since the most
  relevant operators on pseudovarieties do not preserve decidability,
  one often seeks to establish stronger properties. A key role is
  played by relatively free profinite semigroups, which is the
  counterpart of free algebras in universal algebra. The purpose of
  this paper is to give a brief survey of the state of the art,
  highlighting some of the main developments and problems.
\end{abstract}

\maketitle


\section{Why should we be interested in pseudovarieties of semigroups?}
\label{sec:why-pseudovarieties}

In a general algebraic framework, meaning that the operations
considered are fixed, \emph{pseudovarieties} are classes of finite
algebras which are closed under taking homomorphic images, subalgebras
and finite direct products. Although this is a natural analog for
Birkhoff varieties when one is interested in finite algebras, it
lacked external motivation to be investigated until the seminal work
of Eilenberg showing how algebraic structures play a key in the theory
of automata and formal languages \cite{Eilenberg:1974,Eilenberg:1976}.
Earlier results had shown that regular (word) languages over finite
alphabets have certain combinatorial properties if and only if their
syntactic semigroups have corresponding algebraic properties:
star-free versus aperiodic semigroups (meaning that all subgroups are
trivial) \cite{Schutzenberger:1965}, locally testable versus local
semilattices \cite{Brzozowski&Simon:1973,McNaughton:1974}, piecewise
testable versus \Cl J-trivial semigroups \cite{Simon:1975}. Eilenberg
\cite{Eilenberg:1976} showed how all three results fit in a general
bijective correspondence between classes of regular languages with
suitable properties and pseudovarieties of semigroups: the class of
languages corresponding to a given pseudovariety of semigroups \pv V
simply consists of all regular languages over finite alphabets whose
syntactic semigroups belong to~\pv V.

The above three earlier instances of Eilenberg's correspondence are
particularly important because the pseudovarieties in question are
\emph{decidable} in the sense that, for each of them, there is an
algorithm that, given a finite semigroup $S$ (say, by its
multiplication table), determines whether or not $S$ belongs to the
pseudovariety. Since the syntactic semigroup of a regular language
(say, given by a regular expression or a finite automaton recognizing
it) is effectively computable, we see that a pseudovariety is
decidable if and only if so is the corresponding class of regular
languages.

Thus, not only pseudovarieties serve as classifying algebraic
invariants for natural classes of regular languages but they also
translate to algebraic problems combinatorial problems on regular
languages. This provided strong external motivation for the study of
pseudovarieties of semigroups and led to several books dedicated to
the subject
\cite{Lallement:1979,Pin:1986;bk,Almeida:1994a,Rhodes&Steinberg:2009qt}.
Various extensions of Eilenberg's correspondence have been found, such
as: replacing semigroups by arbitrary algebraic structures and
languages by subsets of free algebras in varieties
\cite{Almeida:1990c,Almeida:1994a}, in particular tree languages (sets
of terms) \cite{Gecseg&Steinby:1984,Steinby:1990}; enriching the
algebraic structures with a partial order and/or retaining the
information about the alphabet as generating set for the syntactic
algebra (see \cite{Straubing&Weil:2021}); or even putting the
correspondence in a general categorical setting
\cite{Urbat&Adamek&Chen&Milius:2017,Adamek&Milius&Myers&Urbat:2019}.
In this latter work, duality plays a key role. The importance of Stone
duality had already been pointed out in
\cite[Section~3.6]{Almeida:1994a} and is behind the success of the
profinite approach, but it has only been systematically explored more
recently (see \cite{Gehrke&vanGool:2024}). One important direction
that the research in this area has been exploring is to extend the
Eilenberg theory to deal with classes of languages that are not
necessarily regular. The aim is to separate classes of languages. Such
classes may, for instance, be defined by the complexity of the
membership problem for its members, such as the famous classes
$\mathbf{P}$ and $\mathbf{NP}$ \cite{Fortnow:2019}. Various
recognition devices have been considered, from Stone topological
algebras \cite{Almeida&Klima:2024a}, in particular minimal compact
automata \cite{Steinberg:2013}, to typed
monoids~\cite{Gehrke&Krebs:2017}.

\section{The beginnings}
\label{sec:beginnings}

The remainder of the paper deals mostly with pseudovarieties of
semigroups, although sometimes it is convenient to deal with an
enriched algebraic structure. We start with a brief history of the
subject that predates the motivation arising from the connections with
computer science framed by Eilenberg's correspondence.

Fix an \emph{algebraic signature}, by which we mean a set graded by
the natural numbers $\Sigma=\biguplus_{n\ge0}\Sigma_n$, where the
elements of $\Sigma_n$ are called the \emph{$n$-ary operation
  symbols}. A \emph{$\Sigma$-algebra} is a nonempty set $A$ together
with an \emph{evaluation mapping} $E_n:\Sigma_n\times A^n\to A$ for
each $n\ge0$. Usually, we write $o_A(a_1,\ldots,a_n)$ instead of
$E_n(o,a_1,\ldots,a_n)$.

In case $A$ and each $\Sigma_n$ is further endowed with a topology,
then we say that $A$ is a \emph{topological $\Sigma$-algebra} if each
mapping $E_n$ is continuous.

By a \emph{variety} we mean a class of ($\Sigma$-)algebras closed
under the operators $H$ (adding all homomorphic images), $S$ (adding
all subalgebras), and $P$ (adding all direct products). We say the
algebra $A$ \emph{divides} the algebra $B$, if $A\in HS\{B\}$. An
\emph{identity} is a formal equality $u=v$ between members of an
absolutely free algebra $F_\Sigma X$; we say that it is \emph{valid}
in the algebra $A$ or that $A$ \emph{satisfies} $u=v$ if, for every
homomorphism $\varphi:F_\Sigma X\to A$, the equality
$\varphi(u)=\varphi(v)$ holds. For a set $I$ of identities, we say
that the class $[I]$ consisting of all algebras $A$ that satisfy all
identities from~$I$, written $A\models I$, is \emph{defined by}~$I$;
such classes of algebras are also known as \emph{equational classes}.
In 1935, Birkhoff \cite{Birkhoff:1935} proved that varieties are
precisely the equational classes.

Except for the variety of singleton algebras, varieties contain
infinite algebras and many interesting classes of finite algebras
(such as, in the language with just one binary operation: groups,
nilpotent or solvable groups; nilpotent, completely regular, or
aperiodic semigroups) are not the subclasses of all finite members of
some variety of algebras, that is they are not \emph{equational} as
classes of finite algebras. The natural notion that emerged is that of
a \emph{pseudovariety}, which is a class of finite ($\Sigma$-)algebras
closed under the operators $H$, $S$, and $P_{\mathsf{fin}}$ (adding
all finite direct products).

In 1976, Baldwin and Berman \cite{Baldwin&Berman:1976} showed that, in
case $\Sigma$ is finite, pseudovarieties are the classes of finite
algebras in the union of a countable (increasing) chain of varieties.
In 1985, Ash \cite{Ash:1985} further showed that the finiteness
assumption on $\Sigma$ is superfluous provided chains are replaced by
upper directed unions. Meanwhile, also in 1976, Eilenberg and
Schützenberger \cite{Eilenberg&Schutzenberger:1976} proved that, in
case $\Sigma$ is finite, pseudovarieties are \emph{ultimately defined}
by sequences $(u_n=v_n)_{n\ge0}$ of identities, namely as
\begin{displaymath}
  \bigcup_{n\ge0}\bigcap_{m\ge n}\op u_n=v_n\cl,
\end{displaymath}
where, more generally, for a set $I$ of identities, $\op I\cl$ denotes
the class of all finite algebras that satisfy all identities from~$I$.

In 1981, Reiterman \cite{Reiterman:1982} came up with what proved to
be a rather more fruitful alternative description, which we proceed to
describe. In case $\Sigma$ and $X$ are finite, let $\widehat{F}_\Sigma
X$ be the completion of $F_\Sigma X$ with respect to the following
(ultra)metric
\begin{align*}
  d(u,v)
  &=2^{-r(u,v)}\ \text{where}\\
  r(u,v)
  &=\min\{|A|: A\ \text{finite $\Sigma$-algebra such that }
    A\not\models u=v\}.
\end{align*}
A formal equality $u=v$ with $u,v\in \widehat{F}_\Sigma X$ is called a
\emph{pseudoidentity}. For a finite algebra $A$, $A\models u=v$ means
that, for every homomorphism $\varphi:F_\Sigma X\to A$, the equality
$\hat{\varphi}(u)=\hat{\varphi}(v)$ holds for the unique extension of
$\varphi$ to a continuous mapping $\hat{\varphi}:\widehat{F}_\Sigma
X\to A$, where $A$ is endowed with the discrete topology. Reiterman
proved that 
a class of finite algebras is a pseudovariety if and only if it is
defined by a set $P$ of pseudoidentities:
\begin{displaymath}
  \op P\cl = \{A\ \text{finite $\Sigma$-algebra}: A\models P\}.
\end{displaymath}

Banaschewski \cite{Banaschewski:1983} dropped the finiteness
assumption on~$\Sigma$ by considering instead of a metric structure on
$F_\Sigma X$ a uniform structure, namely that given by taking as a
fundamental system of entourages the set of all congruences of finite
index. The completion considered by Banaschewski may be realized as
the closed subalgebra of the direct product of all quotients $F_\Sigma
X/\theta$ by congruences of finite index consisting of all elements
$(u_\theta)_\theta$ such that, if $\theta_1\subseteq\theta_2$ and
$\varphi:F_\Sigma X/\theta_1\to F_\Sigma X/\theta_2$ is the natural
homomorphism, then $\varphi(u_{\theta_1})=u_{\theta_2}$. The topology
considered in the product is the product topology, where finite
algebras are viewed as discrete topological spaces. Thus, the
completion is an inverse limit of finite algebras, that is, a
\emph{profinite algebra}. Equivalently, it is a compact topological
algebra which is residually finite. For a pseudovariety \pv V, we also
say that a profinite algebra is a \emph{pro-\pv V algebra} if it is an
inverse limit of algebras from~\pv V.

In particular, profinite algebras are topological algebras with
underlying Stone spaces (compact zero-dimensional spaces). Following
\cite{Almeida&Klima:2024a}, we call such algebras $S$ \emph{Stone
  topological algebras}. As recognition devices of subsets of
relatively free algebras $F$, we are interested in the pre-images of
clopen subsets of~$S$ under continuous homomorphisms $F\to S$. The
fact that profinite algebras recognize the same sets as finite
algebras, explains the role played by profinite algebras in finite
semigroup theory. On the other hand, the fact that all Stone
topological semigroups are profinite \cite{Numakura:1957} shows that
one needs to get out of the realm of semigroups to cover non-regular
languages. Exactly when a Stone topological algebra is residually
finite, that is, when it is profinite, is a question that has been
recently solved with various alternative characterizations (see
\cite{Schneider&Zumbragel:2017,Almeida&Klima&Goulet-Ouellet:2023}).
Among such characterizations is the condition that, for every clopen
subset, the syntactic congruence is determined by a finite set of
linear polynomials. Perhaps to understand what this condition means,
it is worth recalling the easy fact that, on an arbitrary monoid $M$,
the linear polynomials $sxt$ determine the syntactic congruence of a
subset $L$ in the sense that the elements $m,n\in M$ are syntactically
equivalent precisely when, for all $s,t\in M$, $smt$ belongs to $L$ if
and only if so does $snt$. Here, the linear polynomial $sxt$ is viewed
as a transformation of $M$ defined by $x\mapsto sxt$, where $s,t\in M$
are fixed and the syntactic congruence $\sigma_L$ of a subset $L$ of
an algebra $A$ is the largest congruence on~$A$ saturating $L$.

As argued in Section~\ref{sec:why-pseudovarieties}, the typical goal
in the application of the theory of pseudovarieties of semigroups in
computer science is to find, if possible, effective characterizations
of the pseudovariety associated by Eilenberg's correspondence to a
suitable class of regular languages. To be more precise, such classes
of languages are called \emph{varieties of languages} and are
correspondences \Cl V associating with each finite alphabet $A$ a
Boolean subalgebra $\Cl V(A)$ of the powerset Boolean algebra $\Cl
P(A^+)$, which consists only of regular languages, such that:
\begin{itemize}
\item if $L\in\Cl V(A)$ and $a\in A$, then $a^{-1}L,
  La^{-1}\in\Cl V(A)$;
\item if $\varphi:A^+\to B^+$ is a homomorphism and $L\in\Cl
  V(B)$ then $\varphi^{-1}(L)\in\Cl V(A)$.
\end{itemize}
The associated pseudovariety \pv V according to Eilenberg's
correspondence is the pseudovariety of semigroups generated by all
syntactic semigroups $A^+/\sigma_L$ of languages $L\in\Cl V(A)$ with
arbitrary finite alphabets $A$. The inverse correspondence can be
described by saying that $\Cl V(A)$ consists of all languages
$L\subseteq A^+$ such that $A^+/\sigma_L$ belongs to~\pv V. If one is
interested in languages that may contain the empty word, the role of
the free semigroup $A^+$ is then played by the free monoid $A^*$, and
instead of semigroups one works with monoids. It is easy to see that a
pseudovariety of monoids consists of the monoids in the pseudovariety of
semigroups it generates.

Now, there are too many pseudovarieties of semigroups for all of them
to be decidable, as there are only countably many algorithms according
to Church's thesis. Indeed the set of all pseudovarieties of
semigroups has at most the power of the continuum as, up to
isomorphism, there are only countably many finite semigroups. But,
even the set of pseudovarieties consisting of finite Abelian groups
has the power of the continuum: for any set $P$ of prime integers, the
cyclic groups of order a member of~$P$ generate a pseudovariety (of
semigroups) for which the set of orders of its cyclic groups of prime
order is precisely $P$; there are a continuum of such sets of primes.

One of the fruitful extensions of the Eilenberg correspondence
involves pseudovarieties of finite ordered semigroups, and positive
varieties of languages \cite{Pin:1995a}. By an \emph{ordered
  semigroup} we man a semigroup with a partial order compatible with
the semigroup operation; homomorphisms between such structures are
required to preserve the operation and the order. \emph{Positive
  varieties of languages} generalize varieties of languages by
replacing Boolean subalgebras of $\Cl P(A^*)$ by 0,1-sublattices.
Pseudovarieties of ordered semigroups may also be viewed as
generalizations of pseudovarieties of semigroups as the mapping
sending the latter to the pseudovariety of ordered semigroups
generated by~\pv V is injective.

\section{Operators on pseudovarieties}
\label{sec:operators}

Several natural algebraic constructions lead to operators on
pseudovarieties. There are basically two kinds of such operators:
\begin{itemize}
\item those that give an explicit, constructive, necessary and
  sufficient condition for a finite semigroup to belong to the
  resulting pseudovariety; and
\item those that define the resulting pseudovariety in terms of
  generators.
\end{itemize}
Some authors prefer to use a different type of letter to indicate
operators but it is most common to use the same type of letter as for
pseudovarieties.

Examples of operators of the first kind:
\begin{itemize}
\item the operator \emph{\pv L}: \pv{LV} is the class of all finite
  semigroups $S$ such that the monoid $eSe$ belongs to~\pv{V} for every
  idempotent $e\in S$;
\item the operator \emph{\pv D}: \pv{DV} is the class of all finite
  semigroups whose regular \Cl D-classes are subsemigroups from \pv V;
\item the operator \emph{\pv E}: \pv{EV} is the class of all finite
  semigroups whose idempotents generate a subsemigroup which belongs
  to~\pv V;
\item the operator \emph{``bar''}: for a pseudovariety of groups \pv
  H, $\overline{\pv H}$~is the class of all finite semigroups whose
  subgroups belong to~\pv H;
\item the operator \emph{\pv B}: \pv{BV} is the class of all finite
  semigroups whose \emph{blocks} are semigroups from~\pv V.
\end{itemize}
For the last example, a \emph{block} of a finite semigroup $S$ is
defined as follows. Consider the smallest equivalence relation $\rho$
on the union $U$ of the subgroups of $S$ that contains the restriction
of both of Green's relations \Cl L and \Cl R to~$U$. A block of $S$ is
obtained by taking the Rees quotient $T/I$ of the subsemigroup $T$
of~$S$ generated by a $\rho$-class by the ideal $I$ consisting of the
elements of $T$ that do not lie in the \Cl D-class of the generators.

Note that all the operators in the preceding paragraph preserve
decidability as all the constructions involved can be effectively
carried out in a given finite semigroup.

Examples of operators of the second kind:
\begin{itemize}
\item the \emph{join} operator: $\pv V\vee\pv W$ is the pseudovariety
  generated by the union $\pv V\cup\pv W$; this consists of all
  divisors of direct products $S\times T$ with $S\in\pv V$ and
  $T\in\pv W$;
\item the \emph{semidirect product} operator: $\pv V*\pv W$ is the
  pseudovariety generated by all semidirect products $S*T$ with
  $S\in\pv V$ and $T\in\pv W$ and consists of all divisors of such
  semidirect products; this is an associative operator
  \cite{Eilenberg:1976};
\item the \emph{Mal'cev product} operator: $\pv V\malcev\pv W$ is the
  pseudovariety generated by the finite semigroups $S$ for which there
  exists a homomorphism $\varphi:S\to T$ such that $T\in\pv W$ and
  $\varphi^{-1}(e)\in\pv V$ for every idempotent $e\in S$;
\item the \emph{power} operator \pv P: \pv{PV} is the pseudovariety
  generated by all semigroups $\Cl P(S)$ of subsets of $S$, under the
  usual subset multiplication $AB=\{ab:a\in A,\ b\in B\}$, with
  $S\in\pv V$;
\item the operator \emph{\pv M}: \pv{MV} is the pseudovariety
  generated by all monoids $S^1$ with $S\in\pv V$, obtained from $S$
  by adding an identity element if $S$ is not a monoid and taking $S$
  itself otherwise; the pseudovarieties in the image of \pv M are
  called \emph{monoidal}.
\end{itemize}
It is known that the first four of these operators do not preserve
decidability
\cite{Albert&Baldinger&Rhodes:1992,Auinger&Steinberg:2001b}. It seems
to be an open problem whether the much less studied operator \pv M
(see~\cite[Chapter~7]{Almeida:1994a}) preserves decidability.

The join, semidirect product, Mal'cev product, and power operators are
all known to correspond to natural operators on varieties of
languages. Thus, it is of interest to be able to relate concrete
values of these operators with decidable operators. There are many
important results of this kind. Such results involve some natural
examples of pseudovarieties which we proceed to list:
\begin{itemize}
\item \pv I: singleton semigroups
\item \pv S: finite semigroups
\item \pv N: finite nilpotent semigroups
\item \pv D: finite semigroups in which idempotents are right
  zeros
\item \pv{Sl}: finite semilattices
\item \pv B: finite bands
\item \pv{Ab}: finite Abelian groups
\item \pv G: finite groups
\item $\pv G_p$: finite $p$-groups
\item $\pv G_{p'}$: finite $p'$-groups (no elements of order $p$)
\item $\pv G_{\mathrm{nil}}$: finite nilpotent groups
\item $\pv G_{\mathrm{sol}}$: finite solvable groups    
\item \pv A: finite aperiodic semigroups
\item \pv J: finite \Cl J-trivial semigroups
\item \pv R: finite \Cl R-trivial semigroups
\item \pv{Com}: finite commutative semigroups
\item \pv{CS}: finite (completely) simple semigroups      
\item \pv{CR}: finite completely regular semigroups
\item $\pv C_n=(\pv{A*G})^n*\pv A$.
\end{itemize}
In the last example, the $n$-th power is taken with respect to the
semidirect product, and $\pv C_0$ is interpreted to be~\pv A. The
pseudovariety $\pv C_n$ is the class of all finite semigroups with
\emph{Krohn-Rhodes complexity} at most $n$
(see~\cite[Chapter~4]{Rhodes&Steinberg:2009qt}). Thus, Krohn-Rhodes
complexity of a finite semigroup $S$ is the minimum $n$ such that
$S\in\pv C_n$; it is computable if and only if all the pseudovarieties
$\pv C_n$ are decidable.

\section{Locality}
\label{sec:locality}

One of the historical instances of Eilenberg's correspondence is that
the pseudovariety corresponding to the variety of all piecewise
testable languages is precisely \pv{LSl}. Here, a language over a
finite alphabet $A$ is said to be \emph{piecewise testable} if
membership of a word $w$ in it can be tested by determining the
factors of length at most some integer $n$ of the word $\$w\#$, where
the symbols $\$$ and $\#$ do not belong to~$A$. The proof that
\pv{LSl} is the right pseudovariety envolves two steps: to show that
the right pseudovariety is \pv{Sl*D} and that $\pv{LSl}=\pv{Sl*D}$.
Both steps have been extensively generalized. For the first step, see
\cite{Straubing:1985}, and~\cite[Chapter~10]{Almeida:1994a}. The
pseudovariety equation $\pv{LV}=\pv{V*D}$ was already considered by
Eilenberg, who called a pseudovariety \pv V \emph{local} if it is a
solution. Examples of local pseudovarieties:
\begin{itemize}
\item \pv R, \pv{R*G} \cite{Stiffler:1973}
\item every nontrivial pseudovariety of groups \cite{Straubing:1985}
\item \pv{CR} \cite{Jones&Szendrei:1992}
\item $\op x^{n+1}=x\cl$ \cite{Jones:1993}
\item \pv{DS} \cite{Jones&Trotter:1995}
\item \pv{DA} \cite{Almeida:1996c}
\item \pv{DG} \cite{Kadourek:2008}
\end{itemize}
The following are examples of non-local pseudovarieties:
\begin{itemize}
\item \pv J \cite{Knast:1983a,Knast:1983b}
\item \pv{Com} \cite{Therien&Weiss:1985}
\item \pv{DAb} \cite{Kadourek:2008}
\item $\pv C_n$ ($n>0$) \cite{Rhodes&Steinberg:2006}
\end{itemize}
In case \pv V is not local, the question remains as to what \pv{V*D}
may be. Straubing \cite{Straubing:1985} showed that, if \pv V is a
decidable monoidal pseudovariety containing the Brandt semigroup $B_2$
(which is the syntactic semigroup of the language $(ab)^+$ over the
alphabet $\{a,b\}$), then \pv{V*D} is also decidable. This restricts
the problem to the subpseudovarieties of~\pv{DS} as this is the
largest pseudovariety that does not contain~$B_2$
\cite{Margolis:1981}. Tilson \cite{Tilson:1987} gave a non-effective
characterization of local (monoidal) pseudovarieties in terms of
\emph{pseudovarieties of categories} which is used in the proof of
some of the above results.

It seems to be an open problem whether locality is decidable for
pseudovarieties given say by a finite basis of computable
pseudoidentities. Probably, the answer is negative, as proofs of
locality/non-locality are in general quite hard.

The interest of locality goes beyond the equation $\pv{V*D}=\pv{LV}$.
For instance: if the pseudovariety \pv V is local and decidable and
the pseudovariety \pv W has a strong form of decidability to be
discussed in Section~\ref{sec:tameness}, then \pv{V*W} is decidable.
This follows from results of \cite{Almeida&Weil:1996} and
\cite{Almeida&Steinberg:2000a}.

\section{Some notable equations}
\label{sec:equations}

It is easy to see that, for a finite group $G$, the blocks of the
power semigroup $\Cl P(G)$ are divisors of~$G$. In particular, the
inclusion $\pv{PH}\subseteq\pv{BH}$ holds for every pseudovariety of
groups \pv H. Margolis and Pin \cite{Margolis&Pin:1985} proved more
precisely that, for such \pv H,
\begin{equation}
  \label{eq:M-P}
  \pv{PH}
  \subseteq\pv{J*H}
  \subseteq\pv{J\malcev H}
  \subseteq\pv{BH}. 
\end{equation}
They also proved that, in the case where $\pv H=\pv G$, the first and
third inclusions are actually equalities. On the other hand, if the
pseudovariety \pv V is local, then the equality
\begin{equation}
  \label{eq:VstarG=VmG}
  \pv{V*G}=\pv{V\malcev G}
\end{equation}
holds, see~\cite[Theorem~3.1]{Henckell&Margolis&Pin&Rhodes:1991}. The equality
(\ref{eq:VstarG=VmG}) also holds for some non-local pseudovarieties
such as~\pv J: this was proved by Henckell and Rhodes
\cite{Henckell&Rhodes:1991} based on results of
Knast~\cite{Knast:1983b} and Ash~\cite{Ash:1991}, the latter is a very
deep and seminal paper that will be mentioned further below. That
completed the proof of the equality $\pv{PG}=\pv{BG}$, thereby showing
that \pv{PG} is decidable; see \cite{Pin:1995} for more on the history
of the proof of this equality.

By a \emph{profinite graph} we mean a topological graph which is an
inverse limit of finite graphs. A profinite graph is \emph{profinitely
  connected} if all its finite continuous homomorphic images are
connected. Let $G$ be a profinite group and $X$ be a subset generating
a dense subgroup. The \emph{profinite Cayley graph} $\Gamma_X(G)$ has
$G$ as vertex space and $G\times X$ as edge space, where $(g,x)$ is an
edge from $g$ to $gx$. In a series of deep papers, Auinger and
Steinberg characterized the pseudovarieties of groups \pv H for which
equality holds in each of the first two inclusions in~(\ref{eq:M-P}):
\begin{itemize}
\item $\pv{PH}=\pv{J*H}$ if and only if \pv H is a \emph{Hall}
  pseudovariety, a property which has many equivalent formulations
  perhaps the simplest of which is that every profinitely connected
  subgraph of the profinite Cayley graph of each finitely generated
  free pro-\pv H group containing the ends of an edge must contain the
  edge itself \cite{Auinger&Steinberg:2003a};
\item $\pv{J*H}=\pv{J\malcev H}$ if and only if \pv H is
  \emph{arboreous}, meaning that the profinite Cayley graphs of all
  finitely generated pro-\pv H groups are \emph{tree-like} in the
  sense that, given any two vertices, there is a unique minimal
  profinitely connected subgraph containing them
  \cite{Auinger&Steinberg:2001a}.
\end{itemize}
The equality $\pv{J\malcev H}=\pv{BH}$ turns out to be exceptional as
it holds only when $\pv H=\pv G$. This follows from
\cite{Higgins&Margolis:1999}, where it is shown that \pv G is the only
pseudovariety of groups \pv H such that $\pv
A\cap\pv{ESl}\subseteq\pv{DA\malcev H}$. This was already observed
in~\cite{Steinberg:1999a}, an explicit proof being given
in~\cite[Corollary~3.5]{Almeida&Klima:2017b}.

The computation of the Mal'cev product \pv{V\malcev G} has deserved a
lot of attention. By a \emph{relational morphism} of monoids we mean a
relation $\mu\subseteq M\times N$ with domain $S$ which is also a
submonoid of~$M\times N$. For a pseudovariety of groups \pv H, the
\emph{\pv H-kernel} of a finite monoid $M$ is the intersection $K_{\pv
  H}(M)$ of all $\mu^{-1}(1)$ when $\mu$ ranges over all relational
morphisms $\mu$ from $M$ to groups in~\pv H. A finite semigroup $S$
belongs to $\pv{V\malcev H}$ if and only if $K_{\pv H}(S^1)\cap S$
belongs to~\pv V.

The Rhodes \emph{Type-II Conjecture} states that $K_{\pv G}(M)$ is the
smallest submonoid of~$M$ containing all idempotents that is closed
under the transformations $m\mapsto amb$ with $aba=a$ or $bab=b$
(known as \emph{weak conjugation}). The conjecture was proved
independently by Ash~\cite{Ash:1991} and Ribes and Zalesski\u\i\
\cite{Ribes&Zalesskii:1993a} (in an equivalent group-theoretic
formulation previously obtained by Pin and Reutenauer
\cite{Pin&Reutenauer:1991}. Algorithms for the \pv H-kernel of finite
monoids are known for other pseudovarieties \pv H such as $\pv G_p$
\cite{Ribes&Zalesskii:1993b}, $\pv G_{\mathrm{nil}}$
\cite{Margolis&Sapir&Weil:1999}, and $\pv{Ab}$ \cite{Delgado:1996a}.
Following the approach of Pin and Reutenauer, all these results
involve the computation of the closure of finitely generated subgroups
in the pro-\pv H topology of a relatively free group (see
\cite{Margolis&Sapir&Weil:1999}). The case of $\pv G_{\mathrm{sol}}$
remains open although there has been some recent progress
\cite{Marion&Silva&Tracey:2024a,Marion&Silva&Tracey:2024b}.

Is is easy to see that $\pv{BG}=\pv{EJ}$ \cite{Margolis&Pin:1984e}.
Thus, the pseudovariety \pv J is a solution of each of the following
equations:
\begin{align}
  \label{eq:V*G=VmG}
  \pv{V*G}&=\pv{V\malcev G}\\
  \label{eq:VmG=EV}
  \pv{V\malcev G}&=\pv{EV}\\
  \label{eq:V*G=EV}
  \pv{V*G}&=\pv{EV}.
\end{align}
No complete solutions are known. That the pseudovariety
$\pv{CS}\cap\overline{\pv{Ab}}$ fails all three equations follows from
the work of Zhang \cite{Zhang:1996} on varieties of completely regular
semigroups. That $\pv{A*G}$ is also a non-solution follows from the
fact that $\pv{(A*G)\malcev G}$ contains semigroups of complexity~2
\cite{Rhodes:1977} (it actually contains semigroups of arbitrarily
high complexity \cite[Corollary~1.2]{Rhodes&Steinberg:2006}). For
\eqref{eq:V*G=EV}, only some important solutions are known, such as
the already mentioned \pv J, but also \pv{Sl} \cite{Ash:1987}, \pv{CR}
(see \cite{Almeida&Escada:1998}), \pv{DA} and \pv{DS} and all
solutions within \pv{LI} \cite{Almeida&Escada:1998}. It is also easy
to see that, if \pv V is a solution of~\eqref{eq:V*G=EV} then so is
every pseudovariety in the interval $[\pv V,\pv{EV}]$. For
\eqref{eq:VmG=EV}, it was claimed without proof by Higgins and
Margolis \cite{Higgins&Margolis:1999} that every pseudovariety in the
interval $[\pv{Sl},\pv{DA}]$ is a solution. Shahzamanian and the
author have announced that there are whole intervals with the power of
the continuum consisting of counterexamples to that claim.

\section{Irreducibility}
\label{sec:decompositions}

Some pseudovarieties are decomposable in terms of others using for
instance one of the binary operators $\vee$, $*$, $\malcev$. In a
vague sense, a decomposition of a pseudovariety \pv V reduces the
study of its members to the study of the members of the
pseudovarieties used in the decomposition which, for simplicity, we
call \emph{factors}. Also, even though the operator used in the decomposition
may not preserve decidability, it may be possible to explore stronger
properties of the factors, cf. Section~\ref{sec:tameness}.

This prompted the investigation of those pseudovarieties that are
indecomposable. For instance, a pseudovariety \pv V is
\begin{itemize}
\item \emph{strictly finite join irreducible} (\emph{sfji}) if
  $\pv V=\pv{U\vee W}$ implies $\pv V=\pv U$ or $\pv V=\pv W$;
\item \emph{finite join irreducible} (\emph{fji}) if $\pv
  V\subseteq\pv{U\vee W}$ implies $\pv V\subseteq\pv U$ or $\pv
  V\subseteq\pv W$.
\end{itemize}
Similar notions may be considered for the operators $*$ and~$\malcev$.
Note that fji implies sfji. The pseudovariety \pv N is sfji but not
fji, see~\cite[Corollary~7.3.30]{Rhodes&Steinberg:2009qt}. The
pseudovariety \pv J is not sfji \cite{Almeida:1991c}.

Margolis, Sapir and Weil \cite{Margolis&Sapir&Weil:1995} showed that,
if \pv H is an extension-closed pseudovariety of groups, then
$\overline{\pv H}$ is irreducible in the stronger sense for $\vee$,
$*$, $\malcev$. Klíma and the author showed that the extension-closure
assumption on the pseudovariety of groups \pv H is superfluous
\cite{Almeida&Klima:2016} and that the pseudovariety \pv A, and for a
nontrivial pseudovariety of groups \pv H, the pseudovarieties
\begin{displaymath}
  \overline{\pv H},\
  \pv{CR}\cap\overline{\pv H},\ 
  \pv{DS}\cap\overline{\pv H},\
  \pv C_n\cap\overline{\pv H}
\end{displaymath}
are all fji even in the lattice of pseudovarieties of ordered
semigroups \cite{Almeida&Klima:2020a}.

The proofs of irreducibility of $\overline{\pv H}$
in~\cite{Margolis&Sapir&Weil:1995} and~\cite{Almeida&Klima:2016}
depend on encoding results that started in \cite{Koryakov:1995}.
Recall that a code is a subset $C$ of a free semigroup $A^+$ such that
the induced homomorphism of free semigroups $C^+\to A^+$ is injective.
The theory of such (variable-length) codes is a well-developed area of
combinatorics on words, see \cite{Berstel&Perrin&Reutenauer:2010}.
Similarly, if \pv V is a pseudovariety, $A$ is finite set, and $C$ is
a subset of the relatively free pro-\pv V semigroup $\Om AV$ generated
by $A$, then we say that $C$ is a \emph{\pv V-code} if the unique
continuous homomorphism $\Om CV\to\Om AV$ extending the inclusion
mapping $C\embedsin\Om AV$ is injective. A closed (even clopen)
subsemigroup of $\Om AV$ need not be a free pro-\pv V semigroup.
However, clopen subgroups of free profinite groups are finitely
generated free profinite groups. In case \pv V satisfies no nontrivial
identities, which is the case of any pseudovariety that contains \pv N
or $\pv G_p$ for some prime~$p$, the subsemigroup of~$\Om AV$
generated by $A$ is a free semigroup on~$A$, which we identify
with~$A^+$. Steinberg and the author \cite{Almeida&Steinberg:2008}
proved that, if \pv H is any extension closed pseudovariety of groups,
then the clopen free pro-$\overline{\pv H}$ subsemigroups of~$\Om
A{}\overline{\pv H}$ are precisely the closures of the subsemigroups
generated by a regular code in~$A^+$. Both to investigate
irreducibility problems and for the role of relatively free profinite
semigroups in symbolic dynamics
\cite{Almeida&ACosta&Kyriakoglou&Perrin:2020b}, it would be worthwhile
to characterize all regular \pv V-codes more generally when the
pseudovariety \pv V satisfies no nontrivial identities.

The idea for the irreducibility results in~\cite{Almeida&Klima:2020a}
can be traced back to the work of Rhodes
(see~\cite[Section~4.6]{Rhodes&Steinberg:2009qt}) considering the
double partial action of a semigroup on a \Cl J-class given by left
and right multiplication. As compact semigroups always have a minimum
ideal, on which both left and right actions are total, the question
that arises is when is it faithful. A weaker property is actually used
in~\cite{Almeida&Klima:2020a} to establish irreducibility results,
namely that both actions are faithful outside the minimum ideal and
the double action is also faithful. This property is established for
free pro-\pv V semigroups on more than one generator for several
pseudovarieties \pv V.

One may also ask when a pseudovariety $\pv V(S)$ generated by a
semigroup $S$ is irreducible for operators such as $\vee$ and
$*$. Note that the stronger irreducibility property holds for
$\pv V(S)$ if and only if whenever $S$ divides a product $T\times U$
of finite semigroups, respectively a semidirect product $T*U$, $S$
must divide a direct power of at least one of the factors $T$ and $U$.
In the case of the join, finite semigroups with this property are said
to be \emph{fji}. Lee, Rhodes, and Steinberg have started a systematic
program of determining (up to isomorphism) all finite fji semigroups.
They have done so for all semigroups of order up to five
\cite{Lee&Rhodes&Steinberg:2019} and for all $\Cl J$-trivial
semigroups of order up to six \cite{Lee&Rhodes&Steinberg:2022}.

Further irreducibility results on pseudovarieties of semigroups may be
found in~\cite[Section~7.3]{Rhodes&Steinberg:2009qt}.

\section{Two key problems}
\label{sec:2-key-problems}

\subsection{Krohn-Rhodes complexity}
\label{sec:KR}

A long outstanding problem in finite semigroup theory is the computation of
the \emph{Krohn-Rhodes complexity} of a finite semigroup. Its
formulation in terms of pseudovarieties was already given at the end
of Section~\ref{sec:operators}: for a given finite semigroup, to
determine the smallest natural number $n$ such that $S\in\pv C_n$.
That there is such an $n$ follows from a theorem of Krohn and Rhodes
\cite{Krohn&Rhodes:1965} showing that every finite semigroup divides a
wreath product of finite simple groups dividing it and as many copies
as needed of a 3-element \Cl L-trivial band monoid. The computation of
the Krohn-Rhodes complexity was started in~\cite{Krohn&Rhodes:1968}
and has motivated many developments in the theory of finite
semigroups. In the book
\cite{Margolis&Rhodes&Schilling:arXiv:2406.18477} one finds many
examples arguing how the Krohn-Rhodes complexity measures complexity
of various natural phenomena.

The computation of the Krohn-Rhodes complexity was given a prominent
role in~\cite{Eilenberg:1976}, thus being already recognized by
Eilenberg as a key problem motivating the theory of pseudovarieties of
semigroups. The books \cite{Arbib:1968} and
\cite{Rhodes&Steinberg:2009qt} give two pictures, separated by four
decades, of the developments aiming toward the solution of the
problem. Recently, a complete, positive, solution of the problem has
been announced by Margolis, Rhodes, and Schilling
\cite{Margolis&Rhodes&Schilling:arXiv:2406.18477}, thereby culminating
almost six decades of deep research in which Rhodes has been the main
driving force. Once confirmed, this is a truly remarkable achievement
and may open the door to further applications of finite semigroup
theory.

\subsection{Dot-depth}
\label{sec:dot-depth}

A second key problem in the theory of pseudovarieties of semigroups is
related with Schützenberger's characterization of star-free languages
\cite{Schutzenberger:1965}. Star-free languages (possibly including
the empty word) are languages that can be expressed in terms of single
letter languages, the language reduced to the empty word, and the
empty language, using only binary union, complementation, and
concatenation. In essence, the Kleene star is replaced by
complementation in the definition of regular expression.

Cohen and Brzozowski \cite{Cohen&Brzozowski:1971} proposed a hierarchy
of star-free languages defined roughly by how many nested levels of
concatenation are needed in a star-free expression for a language. The
levels of the hierarchy define varieties of languages, hence a
corresponding decomposition of the pseudovariety \pv A as the union of
a chain $(\pv B_n)_n$ of pseudovarieties. Brzozowski and Knast
\cite{Brzozowski&Knast:1978} proved that the hierarchy is strict.
Thérien \cite{Therien:1981a} and Straubing \cite{Straubing:1985}
proposed a variant of the hierarchy consisting of monoidal
pseudovarieties which is known as the \emph{Straubing-Thérien
  hierarchy} and which has been the main target of investigation.

On the pseudovariety side, all levels beyond zero of the
(Cohen-)Brzozowski hierarchy $(\pv B_n)_n$ are obtained from the
Straubing-Thérien hierarchy $(\pv V_n)_n$ via the the equation $\pv
V_n*\pv D=\pv B_n$. In view of Simon's \cite{Simon:1975} ($\pv V_1=\pv
J$, whence it is decidable), Knast's \cite{Knast:1983a} ($\pv V_1*\pv
D$ is decidable), and Straubing's \cite{Straubing:1985} already
mentioned results, since $B_2$ belongs to $\pv V_2$, we know that $\pv
V_n$ is decidable if and only if and only if $\pv B_n$ is decidable.

More generally, given a positive variety of languages $\Cl V$, define
$\Pol\Cl V$ to be the class of languages defined by letting $\Pol\Cl
V(A)$ to be the closure under binary union of the set of languages of
the form $L_0a_1L_1\ldots a_nL_n$ with the $a_i\in A$ and the
$L_i\in\Cl V(A)$. Then $\Pol\Cl V$ is a positive variety of languages
\cite{Pin:2013}. Adding complementation in $A^*$ to the closure
operator, one gets a variety of languages $\BPol\Cl V$. The refinement
of the Straubing-Thérien starts with $\Cl V_0(A)=\{\emptyset,A^*\}$
and takes $\Cl V_{n+\nicefrac12}=\Pol\Cl V$ and $\Cl V_{n+1}=\BPol\Cl
V$ for each natural number~$n$.

There is an additional motivation to study the Straubing-Thérien
hierarchy coming from the specification of languages by a
model-theoretic description, where words are viewed as linear
relational models, with a unary predicate $P_a$ for each letter $a$ in
the alphabet such that $P_a(i)$ holds if the letter $a$ is at position
$i$ in the word. Additionally, one considers numerical predicates, for
instance $<$ to say that position $i$ comes before position $j$, and
$\le$ to add the possibility of the two positions coinciding. For
instance, the language $(ab)^+$ is defined over the alphabet $\{a,b\}$
by the following first-order sentence:
\begin{align*}
  &\exists x\,
    \Bigl(
    (P_a(x)\vee P_b(x))
    \wedge
    \bigl(
    (\forall y\, (x\le y))\Rightarrow P_a(x)
    \bigr)
    \wedge
    \bigl(
    (\forall y\, (y\le x))\Rightarrow P_b(x)
    \bigr)
    \Bigr)\\
  &\wedge\forall x,y\,
    \Bigl(
    \bigl(
    x < y
    \wedge 
    \forall z\,
    (z\le x \vee y\le z)
    \bigr)
    \Rightarrow
    \bigl(P_a(x)\Leftrightarrow P_b(y)\bigr)
    \Bigr)
\end{align*}
McNaughton and Papert \cite{McNaughton&Papert:1971} showed that the
languages definable by first order sentences are precisely those whose
minimal automata have aperiodic transition monoids which, combined
with \cite{Schutzenberger:1965} means those are the star-free
languages. Thomas~\cite{Thomas:1982} proved that the dot-depth
hierarchy corresponds to the quantifier-alternation hierarchy in terms
of optimal description of star-free languages.

Various attempts have been made to prove decidability of all
pseudovarieties $\pv V_n$ and there is a long history of partial
results. See \cite{Pin:2017a} for a relatively recent account and
various extensions of the core problem discussed here. In passing, let
us just mention that Pin and Straubing \cite{Pin&Straubing:1985} gave
an interesting algebraic description of~$\pv V_2$ as the pseudovariety
generated by all monoids of upper-triangular Boolean matrices, and
also as~\pv{PJ}. The author~\cite{Almeida:1990e} further showed that
the smallest pseudovariety \pv V such that $\pv{PV}=\pv V_2$ is
generated by the (four-element) syntactic semigroup of the language
$a^*bc^*$ over the alphabet $\{a,b,c\}$. Until recently, these, as
many other attempts to prove decidability of~$\pv V_2$ failed. The
best results so far, due to Place and Zeitoun, give a strong form of
decidability for $\pv V_2$ \cite{Place&Zeitoun:2021} which allows them
to establish decidability for $\pv V_3$ \cite{Place&Zeitoun:2024b}.
The proofs of these results are deep and difficult, which seems to
indicate that, with such an approach, it is going to be very hard to
go further up in the hierarchy as, at least so far, going one step up
it has not been possible to deduce the same strong decidability
property as at the previous step, which precludes a general induction
argument.

Starting with a given variety of languages $\Cl V$ instead of the
trivial variety, one can similarly build on top of it a concatenation
hierarchy whose union is the least concatenation closed variety of
languages containing $\Cl V$. Moreover, if one restricts the products
$L=L_0a_1L_1\ldots a_nL_n$ considered in the definition of the
operator $\Pol$ to be \emph{unambiguous}, in the sense that every word
in $L$ has a unique factorization $w_0a_1w_1\ldots a_nw_n$ with the
$w_i\in L_i$, then one gets corresponding operators $\UPol$ and
$\BUPol$ and unambiguous concatenation hierarchies. The study of such
generalizations has been extensively developed and has been
intricately linked with progress on understanding the original
dot-depth hierarchy. On the pseudovariety side, it is also worth
mentioning that the closure under concatenation of \pv V, meaning the
pseudovariety corresponding to the closure under concatenation of the
variety of languages associated with \pv V is precisely $\pv{A\malcev
  V}$ \cite{Straubing:1979a}, while the unambiguous concatenation
analog is given by $\pv{LI\malcev V}$.
    
\section{Tameness}
\label{sec:tameness}

An approach that has led to many decidability results is to try to
find a basis of pseudoidentities for a given pseudovariety for which
it can be effectively checked whether a given finite semigroup
satisfies it. When the pseudovariety is obtained by applying an
operator to other pseudovarieties, this is sometimes achieved provided
the ``factor'' pseudovarieties satisfy suitable hypotheses. The
general notion of tame pseudovariety was conceived to explore this
idea for the semidirect product.

To explain what is involved in this approach, we need to say a bit
more about relatively fee profinite semigroups first. Let \pv V be a
pseudovariety of semigroups and recall that we denote by \Om AV the
free pro-\pv V semigroup on~$A$, which is endowed with a function
$\iota_{A,\pv V}:A\to\Om AV$ such that, for every function
$\varphi:A\to S$ into a pro-\pv V semigroup $S$, there is a unique
continuous homomorphism $\hat{\varphi}:\Om AV\to S$ such that
$\hat{\varphi}\circ\iota_{A,\pv V}=\varphi$. We may then interpret
each $w\in\Om AV$ as a natural operation on a pro-\pv V semigroup $S$,
namely as the function $w_S:S^A\to S$ given by
$w_S(\varphi)=\hat{\varphi}(w)$. Such operations are known as
\emph{implicit operations}. As \Om AV is pro-\pv V, the
interpretations of $w\in\Om AV$ in semigroups from~$\pv V$ completely
determine~$w$. Among such operations are those defined by elements of
the subsemigroup of some~$\Om AV$ generated by~$A$, which are called
\emph{explcit operations}; in particular, when $|A|\ge2$ and $a$ and
$b$ are distinct letters from~$A$, the explicit operation defined by
$ab$ is just the semigroup multiplication. When $\pv V=\pv S$, a set
of implicit operations containing multiplication is known as an
\emph{implicit signature} and thus gives a structural enrichment of
profinite semigroups. Two commonly considered implicit signatures are:
\begin{itemize}
\item $\omega$, consisting of multiplication and the $\omega$-power, a
  unary operation whose interpretation in a finite semigroup $S$ sends
  each element $s$ to its unique idempotent power $s^\omega$;
\item $\kappa$, consisting of multiplication and the $\omega-1$-power,
  a unary operation whose interpretation in a finite semigroup $S$ sends
  each element $s$ to the inverse of $ss^\omega$ in the maximal
  subgroup containing the idempotent $s^\omega$.
\end{itemize}
In general, for an implicit signature $\sigma$, the subalgebra of~$\Om
AV$ generated by the image of~$\iota_{A,\pv V}$ is the free
$\sigma$-algebra in the variety generated by~$\pv V$ and it is denoted
$\oms AV$. For instance, $\omc AG$ is the free group on~$A$.

One natural and classical problem on a pseudovariety \pv V is whether
the \emph{word problem} for $\oms AV$ is decidable, that is, whether
there is an algorithm that, given two $\sigma$-terms $u$ and $v$,
determines whether or not $\pv V\models u=v$. If the word problem for
$\oms AV$ is decidable for every finite set $A$, then we say that the
\emph{$\sigma$-word problem for~\pv V is decidable}. Some notable
examples of decidability of $\kappa$-word problems are those for the
pseudovarieties $\pv J$ \cite{Almeida:1990b,Almeida:1994a} (which is
intimately related with Simon's characterization of piecewise testable
languages \cite{Simon:1975}), \pv A
\cite{McCammond:1999a,Almeida&Costa&Zeitoun:2015}, \pv R
\cite{Almeida&Zeitoun:2003b}, \pv{DA} \cite{Moura:2009b}, \pv{D(D\vee
  G)} \cite{Borlido:2017a}, \pv{LG}
\cite{Costa&Nogueira&Teixeira:2021}, \pv{A\cap ESl}
\cite{Branco&Costa:2021}, \pv{DAb}
\cite{Almeida&Kufleitner&Wachter:2024}. Several other examples would
be worth considering such as \pv{DS} and \pv{DG} as well as
preservation of decidability of the $\kappa$-word problem under
natural operators, like \pv D, \pv L, \pv{\_*D}, \pv{\_*G},
\pv{LI\malcev\_}, and \pv{A\malcev\_}.

A rather different kind of problem that intervenes in the definition
of tameness has to do with the \pv V-solution of systems of
$\sigma$-equations. Given $\sigma$-identities $u_i=v_i$ ($i\in I$)
over a given finite alphabet $X$, a \emph{\pv V-solution (over the
  alphabet $A$)} is an evaluation of the variables $\varphi:X\to\Om
AS$ such that: $\pv V\models \hat{\varphi}(u_i)=\hat{\varphi}(v_i)$
for every $i\in I$, where $\hat{\varphi}$ is the unique extension
of~$\varphi$ to a $\sigma$-algebra homomorphism $F_\sigma X\to\Om AS$.
Such a solution $\varphi$ is a \emph{$\sigma$-solution} if takes its
values in $\oms AV$. The pseudovariety \pv V is said to be
$\sigma$-reducible for the system $u_i=v_i$ ($i\in I$) if the set of
all $\sigma$-solutions over each finite alphabet $A$ is dense in the
subspace of the product space $(\Om AS)^X$ consisting of all $\pv
V$-solutions. The pseudovariety \pv V is \emph{$\sigma$-reducible} for
a set $E$ of systems of $\sigma$-equations if it is $\sigma$-reducible
for all systems in~$E$. Finally, for a recursively enumerable implicit
signature $\sigma$ consisting of computable implicit operations, a
pseudovariety \pv V is \emph{$\sigma$-tame} for $E$ if the
$\sigma$-word problem for~\pv V is decidable and \pv V is
$\sigma$-reducible for~$E$. In particular, the following terminology
is used instead of saying that \pv V is $\sigma$-tame for~$E$:
\begin{itemize}
\item when $E$ consists of all finite systems of $\sigma$-equations,
  we say that \pv V is \emph{completely $\sigma$-tame};
\item when $E$ consists of all finite systems of equations associated
  with finite directed graphs, in which the variables are the vertices
  and the edges and there is an equation $xy=z$ for each edge
  $x\xrightarrow{y}z$, we say that \pv V is \emph{graph
    $\sigma$-tame};
\item when $E$ consists of the systems $x_1=x_2=\cdots=x_n$ ($n\ge2$),
  we say that \pv V is \emph{pointlike $\sigma$-tame},
\item when $E$ consists of the systems $x_1=x_2=\cdots=x_n=x_n^2$
  ($n\ge2$), we say that \pv V is \emph{idempotent pointlike
    $\sigma$-tame}.
\end{itemize}
More generally, such notions can be considered when inequalities are
taken instead of equations and \pv V is a pseudovariety of ordered
semigroups or monoids.

The interest in tameness comes from the observation, proved
in~\cite{Almeida&Steinberg:2000a,Almeida&Steinberg:2000b} that if \pv
V is $\sigma$-tame for a finite system of equations over a finite set
of variables $X$ and a clopen constraint $K_x\subseteq\Om AS$ is given
(since $K_x$ is the topological closure of a regular language
$L_x\subseteq A^+$, it may be described by finite data) for each $x\in
X$, then it is decidable whether the system has a \pv V-solution over
$A$. This decidability condition had previously been introduced
in~\cite{Almeida:1999b} under the name of \emph{hyperdecidability} as
an approach to proving decidability of semidirect products.

For instance, if the system is reduced to the equation $x=y$, then
being able to determine whether such a solution exists implies that
\pv V is decidable. For, given a finite semigroup $S$ and a continuous
homomorphism $\varphi:\Om AS\to S$, we may take as constraints the
sets $K_x=\varphi^{-1}(s)$ and $K_y=\varphi^{-1}(t)$ for a pair of
distinct points $s,t\in S$. The existence of a \pv V-solution of $x=y$
satisfying those constraints means that there is a pseudoidentity
$u=v$, with $u,v\in\Om AS$, which is satisfied by~\pv V and fails
in~$S$, so that $S\notin\pv V$. On the other hand, the non-existence
of a \pv V-solution satisfying those constraints means that the closed
sets $p_{\pv V}(K_x)$ and $p_{\pv V}(K_y)$ are disjoint, where $p_{\pv
  V}:\Om AS\to\Om AV$ is the unique continuous homomorphism such that
$p_{\pv V}\circ\iota_{A,\pv S}=\iota_{A,\pv V}$; in a Stone space,
this means that there is a clopen set separating the two closed sets,
which in turn means that there is a language $L\subseteq A^+$
recognized by a semigroup from~\pv V that contains the language
$K_x\cap A^+$ and is disjoint from $K_y\cap A^+$. Thus, it is
decidable whether the equation $x=y$ has a \pv V-solution satisfying
given clopen constraints if and only if it is decidable whether given
disjoint regular languages over a finite alphabet $A$ may be separated
by a \pv V-recognizable language. This was first observed
in~\cite{Almeida:1996d}. Note finally that, if the languages
$\varphi^{-1}(r)\cap A^+$ with $r\in S$ may all be pairwise separated
by \pv V-recognizable languages then, as they partition $A^+$, they
are themselves \pv V-recognizable, which entails that $S\in\pv V$. In
particular, if \pv V is tame for the equation $x=y$, then \pv V is
decidable.

What may be considered striking is that difficult decidability
problems may be settled by solving a classical word problem plus
proving a non-algorithmic topological property. Although none of these
ingredients may be easy, achieving them usually means reaching a deep
understanding of the pseudovariety in question, which is why the
combined property is called \emph{tameness}. But, of course, there are
various degrees of tameness depending on for what kind of systems we
are able to prove reducibility and how complicated a signature needs
to be considered. Before proceeding with a survey of tameness results,
we give the main motivation that led to the notion of tameness.

\subsection{Tameness and semidirect product}
\label{sec:semidirect-product}

Exploring Tilson's seminal results using pseudovarieties of categories
to describe semidirect products of pseudovarieties of semigroups or
monoids, namely through his Derived Category Theorem
\cite{Tilson:1987}, Weil and the author \cite{Almeida&Weil:1996}
attempted to describe bases of pseudoidentities for such semidirect
products $\pv{V*W}$. We proceed to describe briefly how such bases are
obtained.

There is an analog of Reiterman's Theorem for pseudovarieties of
categories. Categories are viewed as generalizations of monoids, which
in turn are viewed as categories on a single virtual vertex whose
edges are the elements of the monoid. The role of free profinite
monoids $\Om AS$ on a set $A$ is played by free profinite categories
$\Om\Gamma{Cat}$ on a directed graph $\Gamma$
\cite{Jones:1996,Almeida&Weil:1996}. However, extra care needs to be
taken (which was already present in \cite{Jones:1996} but not in
\cite{Almeida&Weil:1996}) when the graph $\Gamma$ has an infinite
vertex set and in fact it was shown in~\cite{Almeida&ACosta:2007a}
that $\Gamma$ may not generate a dense subcategory of
$\Om\Gamma{Cat}$, which was previously taken for granted in several
papers. In fact, symbolic dynamics is used
in~\cite{Almeida&ACosta:2007a} to show that, starting from the graph
$\Gamma$, it may require an arbitrarily large countable ordinal number of
alternations of taking algebraic generation and topological closure
before $\Om\Gamma{Cat}$ is reached. Yet, pseudovarieties of categories
can be defined by formal equalities $(u=v;\Gamma)$ of elements of free
profinite categories $\Om\Gamma{Cat}$ over finite directed graphs
$\Gamma$ starting and ending at the same vertices
\cite{Jones:1996,Almeida&Weil:1996}. Extra care is needed when
semigroupoids (which are like categories but with no requirement for
local identities) are considered, see Problem~6
in~\cite{Rhodes&Steinberg:2009qt}.

The Basis Theorem (5.3) of~\cite{Almeida&Weil:1996} states that, given
pseudovarieties of monoids $\pv V$ and semigroups $\pv W$, if
$\{(u_i=v_i;\Gamma_i): i\in I\}$ is a basis of pseudoidentities for
the pseudovariety of categories~$g\pv V$ generated by \pv V, then the
set of all the following semigroup pseudoidentities is a basis of
pseudoidentities for $\pv{V*W}$:
$\delta(p)\varepsilon(u_i)=\delta(p)\varepsilon(v_i)$, where $p$ is
the common initial vertex of $u_i$ and $v_i$, $\delta$ is a mapping
from the vertex set of~$\Gamma_i$ to $(\Om AS)^1$, and $\varepsilon$
is a continuous mapping from the edge set of~$\Om{\Gamma_i}{Cat}$ to
$\Om AS$ respecting multiplication such that, for every edge
$q\xrightarrow{x}r$, $\pv W\models \delta(q)\varepsilon(x)=\delta(r)$.

Unfortunately, besides sloppiness in handling graphs with infinite
vertex sets and pseudovarieties of semigroupoids, there is a serious
gap in the proof of a key step (Proposition 3.6) in the proof of the
Basis Theorem where an unjustified exchange of quantifiers is
implicitly made. While no counterexample has ever been produced, it
seems rather unlikely that this key ingredient holds in its full
generality. Thus, the Basis Theorem can for now only be used under one
of the following two extra finiteness assumptions:
\begin{enumerate}
\item\label{item:BTEH-1} \pv W is generated by a finite semigroup;
\item\label{item:BTEH-2} \pv V has finite vertex rank, meaning that
  the pseudovariety of categories $g\pv V$ admits a basis of
  semigroupoid identities over graphs with a bounded number of
  vertices.
\end{enumerate}
Noting that \pv V has vertex rank one if and only if \pv V is local,
local pseudovarieties are specially amenable to this approach. Yet,
there are many pseudovarieties of interest with infinite vertex rank,
for instance,
\begin{itemize}
\item letting $\pv{Com}_{m,\alpha}$ be the pseudovariety of
  commutative monoids satisfying the pseudoidentity
  $x^{m+\alpha}=x^m$, where $m$ is a positive integer and $\alpha$ is
  either a positive integer or $\omega$, when $m\ge2$ there is no
  pseudovariety of finite vertex rank in the (uncountable) interval
  $[\pv{Com}_{m,1},\pv{Com}_{m,\omega}]$ \cite{Almeida&Azevedo:1997};
\item several other intervals of pseudovarieties of infinite vertex
  rank are given in~\cite{Steinberg:2004a} and in the review by
  Auinger of this paper in MathSciNet, MR2025914, including the
  interval between the pseudovariety generated by $B_2^1$ and
  $\pv{DA*H}$ when $\pv H$ is a proper nontrivial pseudovariety of
  groups; among pseudovarieties of infinite vertex rank covered by the
  paper and the review, one finds pseudovarieties such as
  $\pv{\overline{H}*G}$ when $\pv H\subsetneqq\pv G$ (in particular,
  $\pv{A*G}$) and $\pv{Sl*H}$ when $\pv I\ne\pv H\subsetneqq\pv G$.
\end{itemize}

In view of the above discussion, using the tameness approach, the best
that can be stated at present is that if \pv V is a decidable
pseudovariety of monoids of finite vertex rank and \pv W is graph
tame, then $\pv{V*W}$ is decidable \cite{Almeida&Steinberg:2000a}.

\subsection{Tameness and Mal'cev product}
\label{sec:tame-Malcev}

There is also a Basis Theorem for Mal'cev products $\pv{V\malcev W}$
\cite{Pin&Weil:1996a}. The theorem
states that, if the set
$\{u_i(x_1,\ldots,x_{n_i})=v_i(x_1,\ldots,x_{n_i}):i\in I\}$ is a
basis of pseudoidentities for~$\pv V$, then the following is a basis
for $\pv{V\malcev W}$:
\begin{displaymath}
  u_i(w_1,\ldots,w_{n_i})=v_i(w_1,\ldots,w_{n_i})
  \quad \text{whenever } \pv W\models w_1^2=w_1=\cdots= w_{n_i}
  \ (i\in I).
\end{displaymath}
As a corollary, one gets that if $\pv V$ is decidable and $\pv W$ is
idempotent pointlike tame, then $\pv{V\malcev W}$ is decidable.

Thus, the Mal'cev product turns out to be much easier to handle than
the semidirect product by the tameness approach. In contrast, there is
a representation theorem for~$\Om A{(V*W)}$ \cite{Almeida&Weil:1995a}
but no such representation is known for $\Om A{(V\malcev W)}$.

\subsection{Tameness results}
\label{sec:tameness-results}

The following is a summary of known tameness results so far.
\begin{itemize}
\item In seminal work of Ash \cite{Ash:1991} it was proved a property
  that turns out to be equivalent to \pv G being graph $\kappa$-tame
  (see \cite{Almeida&Steinberg:2000a}). Yet, it follows from
  \cite{Coulbois&Khelif:1999} that \pv G is not completely
  $\kappa$-tame, which leads to the question as to whether \pv G is
  completely tame for some signature. As observed
  in~\cite{Almeida&Delgado:1997,Almeida&Delgado:1999}, Ash's result
  turns out to have an interesting formulation in model theory, where
  it was, in that sense, rediscovered by Herwig and Lascar
  \cite{Herwig&Lascar:1997}.

\item The pseudovariety \pv J is completely $\omega$-tame
  \cite{Almeida:2003cshort}. Proving that \pv J is hyperdecidable
  without going through tameness turns out to be much more complicated
  \cite{Almeida&Zeitoun:1998}. But, in fairness, it should be
  mentioned that the algorithm that comes from the tameness approach
  is totally impractical as it involves generating in parallel all
  favorable and unfavorable cases until the one of interest is
  produced.
    
\item If $\pv W$ is such that finitely generated free pro-$\pv W$
  semigroups are finite and computable and $\pv V$ is graph
  hyperdecidable then so is $\pv{V*W}$ \cite{Almeida&Silva:1997a}.
  This should be improvable to tameness but does not appear to have
  been done so far.

\item The pseudovariety \pv{CR} is $\kappa$-tame for graph systems of
  equations. This follows from~\cite{Almeida&Trotter:2001} together an
  observation of K. Auinger that the required supposedly improved
  tameness of \pv G is actually granted by Ash's result. There is also
  a potential problem with the proof because of the usage of free
  profinite categories over infinite-vertex graphs, which are assumed
  to generate dense subcategories. Yet, it was shown in
  \cite{Almeida&ACosta:2015a} that the required property does hold in
  the case in question.

\item The pseudovariety $\pv G_p$ is not graph $\kappa$-tame but it is
  graph $\sigma$-tame for a certain infinite implicit signature
  $\sigma$ constructed using ideas from symbolic dynamics
  \cite{Almeida:1999c}. This depends on results of
  Steinberg~\cite{Steinberg:1998a}, who previously proved a weak form
  of hyperdecidability for~$\pv G_p$.

\item The pseudovariety \pv{Ab} is completely $\kappa$-tame
  \cite{Almeida&Delgado:2001}.
    
\item The pseudovariety \pv R is completely $\omega$-tame
  \cite{Almeida&Costa&Zeitoun:2005b}. The idea of the proof of
  complete reducibility is to adapt that of Makanin's algorithm to
  solve equations in free semigroups \cite{Makanin:1977,Makanin:1981},
  even though in our case there is no algorithm involved.

\item The following was established
  in~\cite{Delgado&Masuda&Steinberg:2007}:
  \begin{itemize}
  \item A monoidal pseudovariety of commutative
    semigroups is completely $\kappa$-hyperdecidable if and only if
    it is decidable.
  \item If a proper pseudovariety of Abelian groups is
    $\kappa$-reducible for systems of graph equations, then it is
    locally finite.
  \end{itemize}
  
\item The pseudovariety \pv{LSl} is completely $\omega$-tame
  \cite{Costa&Nogueira:2009}.

\item The pseudovarieties \pv A and \pv{DA} are pointlike
  $\omega$-tame \cite{Almeida&Costa&Zeitoun:2015b}. In the case of \pv
  A, the proof uses ideas of Henckell
  \cite{Henckell:1988,Henckell&Rhodes&Steinberg:2010} giving an
  algorithm to determine the semigroup of all \pv A-pointlike subsets
  of a finite semigroup.
  
\item The pseudovariety \pv{LG} is graph $\kappa$-tame
  \cite{Costa&Nogueira&Teixeira:2021}. The proof uses the more general
  result that, if \pv V is graph $\kappa$-reducible, then so is
  \pv{V*D} \cite{Costa&Nogueira&Teixeira:2016}.

\item The pseudovariety of groups \pv H is completely $\kappa$-tame if
  and only if so is $\pv{D(D\vee H)}$ \cite{Almeida&Borlido:2017}. The
  method is similar to that used for \pv R, described above, but is
  more complicated because groups are involved.

\item The pseudovariety $\pv G_{\mathsf{nil}}$ is graph tame
  \cite{Alibabaei:2017c}.
    
\item That the pseudovariety \pv{DAb} is completely $\kappa$-tame has
  recently been announced by Kufleitner, Wächter and the author but at
  the moment only a preprint is available for the $\kappa$-word
  problem \cite{Almeida&Kufleitner&Wachter:2024}.
\end{itemize}

The tameness approach has also been explored to compute joins. For
instance, $\pv{J\vee G}$ was independently shown by Steinberg
\cite{Steinberg:1997a} and Azevedo, Zeitoun and the author
\cite{Almeida&Azevedo&Zeitoun:1997}. The former work fits in a more
comprehensive approach to joins, giving hyperdecidability results, and
was part of the author's Ph.D. thesis \cite{Steinberg:1998}. There are
many other papers dealing with joins of pseudovarieties, often using
the profinite, but not necessarily tameness, aproach.
\cite{Almeida:1988b,Almeida&Azevedo:1989,Almeida&Costa&Zeitoun:2004,
  Almeida&Weil:1994c,Costa:2002a,Auinger:2002,Costa:2004,Steinberg:2001b}

It is also worth mentioning that the tameness approach has been
proposed to establish decidability of the Straubing-Thérien hierarchy.
This was started by Klíma, Kunc, and the author
\cite{Almeida&Klima&Kunc:2016} by showing that the $\omega$-inequality
problem, meaning solving the inequality $x\le y$, is decidable over
all pseudovarieties $\pv V_n$ and $\pv V_{n+\nicefrac12}$ ($n$
non-negative integer) in the refined hierarchy. This reduces the
decidability of the hierarchy to establishing the purely topological
property of $\omega$-reducibility for the inequality $x\le y$ of all
levels. As evidence that such a property may hold, Vola\v ríková
\cite{Volarikova:2024} has shown that $\pv V_2$ is defined by
$\omega$-identities: if the topological property holds then all levels
of the hierarchy would be defined by $\omega$-inequalities.

\section{The structure of relatively free profinite semigroups}
\label{sec:rel-free}

In view of the role of relatively free profinite semigroups in the
profinite approach, it is worth understanding the structure of such
semigroups. This is in general quite hard and has only been achieved
in very few cases. An idea that has been extensively explored is that,
just as the positions of the letters in finite words are linearly
ordered, members of relatively free profinite semigroups, sometimes
called \emph{profinite words}, but which the author prefers to call
\emph{pseudowords}, should also have some kind of linear structure.
Even for the pseudovarieties of groups this is in a sense the case:
the profinite Cayley graph of $\Om AH$ is a profinite $\pv H$-tree if
and only if $\pv{(H\cap Ab)*H}=\pv H$ \cite{Almeida&Weil:1994c}, which
extends results of Gildenhuys and Ribes \cite{Gildenhuys&Ribes:1978}.

The aperiodic case also presents a linear behavior. This had already
been observed for \pv J \cite{Almeida:1990b}, \pv R
\cite{Almeida&Weil:1996b}, \pv{DA} \cite{Almeida:1996c}, $\pv{D(D\vee
  H)}$ \cite{Almeida&Borlido:2017}, but the order types become much
more involved for the pseudovariety \pv A
\cite{Almeida&Costa&Zeitoun:2009a,Almeida&ACosta&Costa&Zeitoun:2019,Gool&Steinberg:2019}.
The latter of these works also brings about interesting connections
with model theory, which we have seen in this survey to pop up every
so often.

The local structure of relatively free profinite semigroups has also
been investigated, particularly, the structure of their regular $\Cl
D$-classes. The author came up with an interesting connection with
symbolic dynamics: for every pseudovariety \pv V containing
$\pv{LSl}$, the regular $\Cl J$-classes of $\Om AV$ that are maximal
in the partial order of $\Cl J$-classes are in bijection with the
minimal shift spaces $X\subseteq A^{\mathbb{Z}}$ \cite{Almeida:2005c}:
the $\Cl J$-class associated with $X$ consists of all non-finite
pseudowords in $\overline{L(X)}$, where $L(X)$ consists of all finite
words that appear as blocks in elements of~$X$. Recall that a
\emph{shift space} over a finite alphabet $A$ is simply a nonempty
closed subset of $A^{\mathbb{Z}}$, whose elements are viewed as
biinfinite words, which is stable under shifting the origin. We say
that $X$ is \emph{sofic} if $L(X)$ is a regular language,
\emph{irreducible} if, for all $u,v\in L(X)$ there is $w$ such that
$uwv\in L(X)$, \emph{periodic} if $L(X)$ consists of all factors of
the powers of a fixed word, and \emph{substitutive} if $L(X)$ consists
of all factors of $\varphi^n(a)$ where $a\in A$ and $\varphi$ is a
primitive endomorphism of~$A^+$.

More generally, irreducible shift spaces also have a unique $\Cl
J$-minimal $\Cl J$-class intersecting $\overline{L(X)}$, which is
denoted $J_{\pv V}(X)$. Since all maximal subgroups in a $\Cl J$-class
are isomorphic this led to the definition of the \emph{Schützenberger
  group} of~$X$, denoted $G_{\pv V}(X)$, to be any of the maximal
subgroups of $J_{\pv V}(X)$. In the case of a minimal shift space,
there is a natural geometric interpretation of $G_{\pv V}(X)$ as an
inverse limit of profinite completions of Poincaré groups of certain
Rauzy graphs of $X$ \cite{Almeida&ACosta:2016b}. It is also an
invariant of topological conjugacy, which is the natural notion of
isomorphism between shift spaces \cite{Costa:2006} (see also
\cite{ACosta&Steinberg:2016}). The book
\cite{Almeida&ACosta&Kyriakoglou&Perrin:2020b} gives an introduction
to this theory and interesting connections with coding theory which
were already explored
in~\cite{Almeida&ACosta&Kyriakoglou&Perrin:2020,Goulet-Ouellet:2022a}.

A remarkable result of Costa and Steinberg
\cite{ACosta&Steinberg:2011} shows that, whenever \pv H is an
extension-closed pseudovariety of groups and $X$ is an irreducible
sofic shift space, then $G_{\overline{\pv H}}(X)$ is a free
pro-$\overline{\pv H}$ group which is of countable rank unless $X$ is
periodic, in which case the group is procyclic
\cite{Almeida&Volkov:2006}. In contrast, it had already been observed
in~\cite{Almeida:2005c} that $G_{\pv S}(X)$ may not be a free
profinite group even for substitutive shift spaces. In the case of
substitutive shift spaces, a finite (profinite) presentation can be
computed for $G_{\pv S}(X)$ which entails that it is decidable whether
a finite group is a continuous quotient of $G_{\pv S}(X)$ and allows
to prove freeness or non-freeness (even relatively to any
pseudovariety, as is the case for the much studied Prouhet-Thue-Morse
shift space, which is generated by the substitution $a\mapsto ab$,
$b\mapsto ba$) in many cases \cite{Almeida&ACosta:2013}. Further
relevant results for freeness of Schützenberger groups have also been
obtained in~\cite{Goulet-Ouellet:2022d,Goulet-Ouellet:2022c}.

It remains an open problem what kind of profinite group can $G_{\pv
  V}(X)$ be when $X$ is a minimal shift space. More information
provided by $\Cl J$-classes associated with a shift space has also
been explored in~\cite{Almeida&ACosta:2012,ACosta&Steinberg:2021}.

\section{Conclusion}
\label{sec:conclusio
n}

There are many aspects of the theory of pseudovarieties, which extends
for over six decades, that it is impossible to cover in such a brief
survey. By no means this is meant to belittle such aspects and the
many valuable contributions that many authors have made, but rather
reflects the limitations of the author of this survey.

In any case, it is hoped that this work gives a feeling for the
richness and depthness of a well-motivated theory whose potential
applications have perhaps not yet been fully explored.

\section*{Acknowledgments}

The author acknowledges partial support by CMUP (Centro de Matemática
da Universidade do Porto), member of LASI (Intelligent Systems
Associate Laboratory), which is financed by Portuguese funds through
FCT (Fundação para a Ciência e a Tecnologia, I. P.) under the projects
UIDB/00144/2020 and UIDP/00144/2020.

%
%
\bibliographystyle{amsplain}
\bibliography{sgpabb,pseudovarieties-prep}
\end{document}